\documentclass[letterpaper, 10 pt, conference]{ieeeconf}  

\IEEEoverridecommandlockouts                              

\usepackage{amsfonts, dsfont, mathtools, amsmath, amssymb, bbold, mathrsfs, bm, upgreek}
\usepackage{graphicx, float, epsfig, color, psfrag, adjustbox, enumerate}
\usepackage{subcaption}
\usepackage{refcount, url, cleveref}
\usepackage[noadjust]{cite}
\usepackage[usenames,dvipsnames,svgnames,table]{xcolor}
\usepackage{algorithm}
\usepackage{algpseudocode}

\algrenewcommand\alglinenumber[1]{#1:}
\usepackage{multirow}
\usepackage{multicol}
\usepackage{tikz}
\usetikzlibrary{positioning,arrows}

\title{\LARGE \bf
Optimal Control of a Soft CyberOctopus Arm
}

\author{Tixian Wang$^{1,2}$, Udit Halder$^2$, Heng-Sheng Chang$^{1,2}$, 
 Mattia Gazzola$^{1,3,4}$, Prashant G. Mehta$^{1,2}$
\thanks{We gratefully acknowledge financial support from ONR MURI N00014-19-1-2373, NSF/USDA $\#$2019-67021-28989, and NSF EFRI C3 SoRo $\#$1830881. We also acknowledge computing resources provided  by the Blue Waters project (OCI- 0725070, ACI-1238993), a joint effort of the University of Illinois at Urbana-Champaign and its National Center for Supercomputing Applications, and the Extreme Science and Engineering Discovery Environment (XSEDE) Stampede2 (ACI-1548562) at the Texas Advanced Computing Center (TACC) through allocation TG-MCB190004.}%
\thanks{$^{1}$Department of Mechanical Science and Engineering, $^{2}$Coordinated Science Laboratory,
$^{3}$Department of Molecular and Integrative Physiology, $^{4}$National Center for Supercomputing Applications, \&  $^{7}$Carl
R. Woese Institute for Genomic Biology, University of Illinois at Urbana-Champaign. 
  Corresponding e-mail:  {\tt\small udit@illinois.edu}}%
  \thanks{The first author is thankful to Arman Tekinalp for helpful discussions on numerical solvers.}
}
\def\R{{\mathds{R}}}

\def\0{{\mathbb{0}}}
\def\1{{\mathds{1}}}

\def\a{{\mathbf{a}}}
\def\b{{\mathbf{b}}}

\newcommand{\norm}[1]{\left\lVert#1\right\rVert}
\newcommand{\inner}[2]{\left\langle #1, #2 \right\rangle}

\definecolor{db}{RGB}{23,20,119}
\definecolor{dg}{RGB}{2,101,15}

\newtheorem{proposition}{Proposition}[section]

\newtheorem{remark}{Remark}

\usepackage{upgreek}
\newcommand{\dif}{\mathrm{d}}

\newcommand{\set}[1]{\left\{#1\right\}}

\newcommand{\material}[1]{
	\ifthenelse{\equal{#1}{\kappa}}{\upkappa}{
	\ifthenelse{\equal{#1}{\nu}}{\upnu}{
	\ifthenelse{\equal{#1}{\omega}}{\upomega}{
	\ifthenelse{\equal{#1}{\sigma}}{\upsigma}{
	\ifthenelse{\equal{#1}{\theta}}{\uptheta}{
	\mathsf{#1}}}}}}
}

\newcommand{\target}{\text{target}}
\newcommand{\intrinsic}{\circ}

\newcommand{\deformations}{w}
\newcommand{\states}{q}
\newcommand{\momentums}{p}
\newcommand{\costates}{\xi}
\newcommand{\Hamiltonian}{\mathcal{H}}
\newcommand{\Lagrangian}{\mathcal{L}}
\newcommand{\potential}{\mathcal{V}}
\newcommand{\kinetic}{\mathcal{T}}

\renewcommand{\a}{\mathsf{a}}
\renewcommand{\b}{\mathsf{b}}

\newcommand{\Gmuscle}{\mathsf{G}}

\newcommand{\sobolev}{\mathrm{H}}
\newcommand{\lspace}{\mathrm{L}^2}
\newcommand{\controlspace}{\mathrm{U}}
\newcommand{\statespace}{\mathrm{Z}}

\newcommand{\ud}{\,\mathrm{d}}
\newcommand{\diag}{\text{diag}}
\newcommand{\x}{b}
\newcommand{\y}{a}
\newcommand{\pt}{p^r}
\newcommand{\pr}{p^\theta}
\newcommand{\mut}{\mu^r}
\newcommand{\mur}{\mu^\theta}
\newcommand{\gammat}{\gamma^r}
\newcommand{\gammar}{\gamma^\theta}
\newcommand{\uf}{u^F}
\newcommand{\uc}{u^C}
\newcommand{\D}{\tilde{\mathcal{D}}}
\newcommand{\diff}{\bar{\mathcal{D}}}

\begin{document}
\bstctlcite{BSTcontrol} 
\maketitle
\thispagestyle{empty}
\pagestyle{empty}


\begin{abstract}
In this paper, we use the optimal control methodology to control a flexible, elastic Cosserat rod.  An inspiration comes from stereotypical movement patterns in octopus arms, which are observed in a variety of manipulation tasks, such as reaching or fetching. To help uncover the mechanisms underlying these observed morphologies, we outline an optimal control-based framework. A single octopus arm is modeled as a Hamiltonian control system, where the continuum mechanics of the arm is modeled after the Cosserat rod theory, and internal, distributed muscle forces and couples are considered as controls. First order necessary optimality conditions are derived for an optimal control problem formulated for this infinite dimensional system. Solutions to this problem are obtained numerically by an iterative forward-backward algorithm. The state and adjoint equations are solved in a dynamic simulation environment, setting the stage for studying a broader class of optimal control problems. Trajectories that minimize control effort are demonstrated and qualitatively compared with experimentally observed behaviors.  
\end{abstract}

\begin{keywords}
Cosserat rod, optimal control, maximum principle, soft robotics, octopus, Hamiltonian systems
\end{keywords}

\section{Introduction} \label{sec:intro}

\subsection{Background and Objectives}
Over the past few decades, the optimal control paradigm has been increasingly used to explain and understand dynamic phenomena in biological systems. 
Examples range from game theoretic models of population dynamics \cite{smith1973logic, smith1982evolution} to testing optimality hypotheses for collective motion in starling murmurations \cite{attanasi2014information, justh2015optimality, halder2019optimality}, or the minimum-jerk hypothesis for movement planning \cite{flash1985coordination, viviani1995minimum, todorov1998smoothness, todorov2004optimality}. Through a mixture of experimental data analysis and theoretical modeling, these approaches often reveal deep insights into the underlying mechanisms at play \cite{todorov2002optimal, todorov2004optimality}. In this work, we take a similar route to examine the problem of octopus arm movement.

Flexible octopus arms are excellent candidates for studying the intricate interplay between continuum mechanics and sensorimotor control. As opposed to articulated limbs in humans, octopus arms are soft and possess a complex muscular architecture that provides exquisite manipulation control.
The muscles are independently innervated by motor neurons along the arm enabling a rich repertoire of deformations -- stretch, shear, bend, and twist. However, despite their virtually infinite degrees of freedom -- and thus having many options to carry out a single task -- octopuses are observed to engage in certain (task-specific) stereotypical movement strategies.  In experimental studies~\cite{gutfreund1996organization, gutfreund1998patterns, sumbre2006octopuses}, these strategies are broadly categorized into two groups.

\medskip

\begin{figure*}[t]
	\centering
	\includegraphics[width=2.\columnwidth,trim = {0pt 0pt 0pt 0pt}]{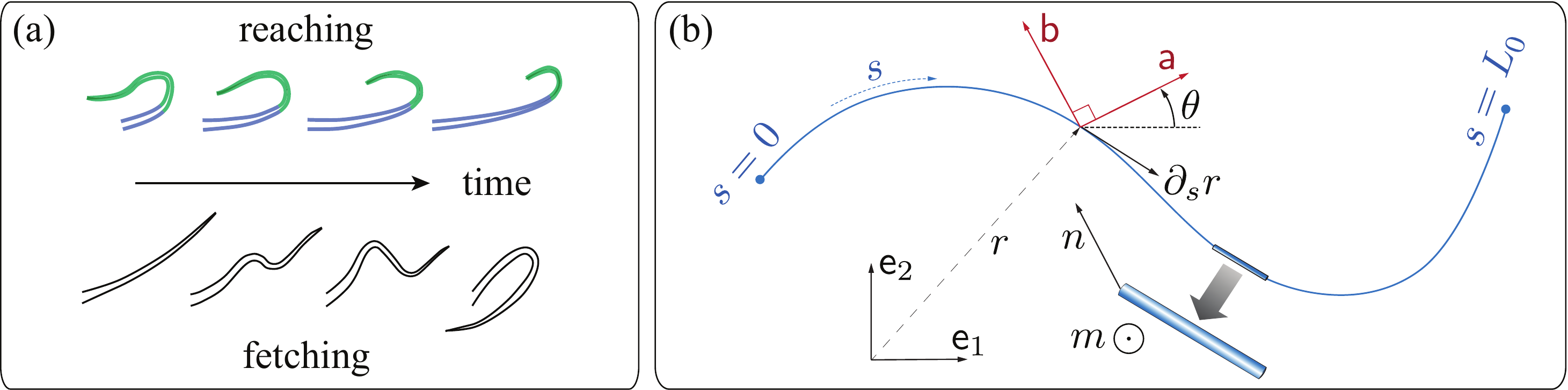}
	\caption{(a) The octopuses have been observed to exhibit bend propagation (for reaching) and elbow formation (for fetching). The bend propagation is actively achieved by propagating muscle actuation, illustrated by blue color; green represents the unactuated portion of the arm. (b) A schematic of the planar Cosserat rod model.}
	\label{fig:rod} 
	\vspace{-10pt}
\end{figure*}

\noindent\textbf{(i) Reaching pattern -- bend propagation:}
For the task of reaching to a fixed target (Fig.~\ref{fig:rod}a), the arm creates a bend at the base of the arm and propagates that bend toward the tip \cite{gutfreund1996organization}. It was later showed that these waves are not mere whip-like mechanical waves \cite{hines1978bend, coleman1992flexure} due to the flexible arm structure, rather the bend propagation is achieved by actively creating  waves of muscle activation signals \cite{gutfreund1998patterns}. Electromyogram (EMG) recordings of muscle activation reveals association of muscle contraction with the traveling bend. Ex-vivo experiments seem to suggest that these movement patterns may actually be encoded in the neural circuitry of the arm itself \cite{sumbre2001control}. 

\medskip

\noindent\textbf{(ii) Fetching pattern -- creation of pseudo-joints:}
The octopus typically employs a different strategy for the scenario of fetching food to its mouth. In this case, the arm behaves like an articulated limb \cite{sumbre2005motor, sumbre2006octopuses} (see Fig.~\ref{fig:rod}a), creating dynamic pseudo-joints at three locations along the arm -- proximal, medial, and distal. The medial joint is the most prominent one, and forms at the location where two waves of propagating muscle activation collide.

\medskip

The objective of the present paper is to introduce an optimal control framework, associated numerical algorithms, and software tools to systematically investigate potential optimality bases of these stereotypical movement strategies.  We are particularly interested in understanding the traveling wave phenomena observed in experimental studies.  The framework introduced here is seen as a first step towards an inverse optimality analysis of the observed behaviors.

\subsection{Contributions}
The dynamics of a soft arm are modeled using the Cosserat rod theory \cite{antman1995nonlinear, gazzola2018forward, chang2020energy}.  Internal muscle forces and couples, when considered as control inputs, give rise to a control system in an infinite-dimensional state space setting. Since the observed stereotypical arm movements occur primarily in-plane \cite{gutfreund1996organization}, we restrict our modeling to planar settings, leading to a control system described by six nonlinear PDEs. We propose an optimal control problem associated with this control system. The Pontryagin's Maximum Principle (PMP) is used to derive the six adjoint PDEs for the costate variables. The PMP is also used to obtain the (open-loop) optimal control input. 

The resulting two-point boundary value problem is numerically solved in an iterative manner, referred to here as the forward-backward algorithm. The forward path, or the Cosserat dynamical equations are solved using the existing software tool~\textit{Elastica} \cite{gazzola2018forward, zhang2019modeling,naughton2021elastica}. A custom solver is implemented to simulate the backward path or the costate equations. The deviation from optimality is utilized to adjust the control in an iterative manner so as to achieve optimality.

The numerical solver is applied to three test cases related to the reaching and the fetching movement patterns. Simulation results are used to qualitatively compare with observed wave propagations or elbow forming.


\subsection{Paper Outline}

The remainder of this paper is organized as follows: In Sec.~\ref{sec:problem}, the Cosserat rod model and dynamics in the planar case are introduced and an optimal control problem is formulated. The solution to the optimal control problem, including the forward-backward algorithm and the numerical methods are described in Sec.~\ref{sec:solution}. Results of numerical experiments appear in Sec.~\ref{sec:numerics}. The paper is concluded in Sec.~\ref{sec:conclusion}.

\section{Problem formulation} \label{sec:problem}

\subsection{Dynamic modeling of an arm as a Cosserat rod}
Let $\set{\mathsf{e}_1,\mathsf{e}_2}$ denote a fixed orthonormal basis for the two-dimensional {laboratory frame}. Time $t\in \R$ and arc-length $s\in[0,L_0]$, $L_0$ being the length of the undeformed rod, represent the two independent variables. The partial derivatives with respect to $t$ and $s$ will be denoted by the subscripts $(\cdot)_t$ and $(\cdot)_s$, respectively. 

The \textit{state} of the rod is described by the vector-valued function $\states (t,s) = (r(t,s) , \theta(t,s))$  
where $r=(x,y)\in\R^2$ denotes the position vector of the centerline, and the angle $\theta\in\R$ defines 
{the material frame} spanned by the orthonormal pairs $\{\a,\b\}$, where $\a = \cos \theta \,\mathsf{e}_1 + \sin \theta \, \mathsf{e}_2, ~ \b = -\sin \theta \, \mathsf{e}_1 + \cos \theta \, \mathsf{e}_2$ (see Fig.~\ref{fig:rod}b).
The vector $\a$ is normal to the cross section.  The deformations $\deformations=(\nu_1,\nu_2,\kappa)$, {stretch, shear, and curvature,} are related to the local frame $\{\a, \b\}$ through $r_s =\nu_1\a+\nu_2\b$ and $\theta_s = \kappa$.  Finally,  $p(t,s)=\mathcal{M}q_t(t,s)$ is used to denote the momentum variable where $\mathcal{M}$ is the mass-inertia density matrix. 

The Hamiltonian formulation requires specification of the kinetic energy $\mathcal{T}$ and the potential energy $\mathcal{V}$ of the rod as follows:
\begin{equation*}
	\mathcal{T}(p) = \frac{1}{2}\int_0^{L_0}  p^{\mathsf{T}} \mathcal{M}^{-1}p ~\dif s,\quad
	\mathcal{V}(q) = \int_0^{L_0}W(w) \ud s
\end{equation*}
where $W:w \mapsto \R$ is referred to as the stored energy function of the rod. A quadratic stored energy function, which leads to a linear stress-strain relationship, is used in this work. The total energy function or the Hamiltonian $\mathcal{H}(q,p) := \kinetic(p)+\potential(q)$ yields the Hamilton's equations of the rod dynamics in the classical Cosserat theory \cite{antman1995nonlinear, chang2020energy}. 

The \textit{generalized state} of the rod is denoted as 
\[
z(t) := (\states(t,\cdot), \momentums(t,\cdot)) \in \statespace, ~ t\in [0, T]
\]  
An appropriate choice of function space is $\statespace=\sobolev^1 ([0, L_0]; \R^3)\times \lspace ([0, L_0]; \R^3)$  equipped with the appropriate boundary conditions. The dynamics of the Hamiltonian control system are expressed as follows:
\begin{align}
\frac{\dif z}{\dif t} (t) = (J - R) \frac{\delta \Hamiltonian}{\delta z} + \Gmuscle(z(t)) u(t) =  : f(z(t), u(t))
\label{eq:control_system_abstract}
\end{align}
where $z(0)$ is the initial condition, $J$ is the skew-symmetric structure matrix $\begin{psmallmatrix*}[c] 0 & \1 \\ -\1 & 0 \end{psmallmatrix*}$, and $R = \begin{psmallmatrix*}[c] 0 & 0 \\ 0 & \zeta \1 \end{psmallmatrix*}$ is the dissipation matrix, $\zeta>0$ is a damping coefficient, modeling viscoelastic effects in the rod \cite{gazzola2018forward}.  The term 
$\Gmuscle(z(t)) u(t)$ on the right hand side is used to model the effect of the distributed internal muscle forces and couples. The functions $u(\cdot) \in \mathfrak{U}$ are called control inputs. Here $\mathfrak{U}$ is the set of all measurable functions $u(\cdot): [0, T] \rightarrow \controlspace$, where $\controlspace$ is a suitable function space called the control space. We take this as the $\lspace ([0, L_0]; \R^3)$ space.
The modeling of $\Gmuscle$ is complicated and depends on the muscle type details of the octopus.  In this paper, we make the simplifying assumption $\Gmuscle (z(t)) \equiv \begin{psmallmatrix} 0 \\ \1 \end{psmallmatrix}$. 

The explicit form of the six partial differential equations in the model \eqref{eq:control_system_abstract} appears in Appendix \ref{appdx:explicit}.

\subsection{An optimal control problem}

Both sterotypical movement patterns introduced in Sec.~\ref{sec:intro} involve reaching a given target point $\states^\target \in \R^3$.  Even if realistic muscle constraints were considered (they are ignored here), there would exist a large number of potential strategies to achieve the objective.  Optimal control appears to be a natural choice to obtain a unique strategy.  This is done through      
formulating the
following free endpoint optimal control problem:
\begin{align}
	\begin{split}
		\underset{u}{\text{minimize}} ~ \mathcal{J}({u}) &= \int_0^T  \Lagrangian(z(t), u(t)) \ud t + \Phi (z(T))\\ 
		\text{subject to} ~ &\text{\eqref{eq:control_system_abstract} and a given}~ z(0,s)
	\end{split}
		\label{eq:optimal_control_problem_free_end}
\end{align}
Here the end point $z(T) = (\states(T), \momentums(T))$ is free and penalizes the cost $\Phi$ associated with the underlying task, for example the distance from the arm tip to the designated target point. Note that a free endpoint problem is considered as opposed to a fixed endpoint problem due to the ease in algorithmic implementation as described in Sec. \ref{sec:forward_backward}.

The choice of the cost function is problem dependent.   In this paper, a quadratic model is assumed for the control cost and the elastic potential energy is assumed for the state-dependent cost
\begin{align}
\Lagrangian(z, u) = \frac{1}{2} \| u\|^2_{\lspace} + \chi_1 \potential(\states)
\label{eq:Langrangian}
\end{align}
where the weighting parameter $\chi_1>0$ is used to penalize the deformation of the arm.  The terminal cost is used in place of a fixed endpoint constraint
\begin{align}
\Phi(z(T)) =  \chi_2 \Phi_{\text{tip}}({q}(T, L_0), \states^\target)
\label{eq:terminal_cost}
\end{align}
where the function $\Phi_{\text{tip}}$ measures the distance between the arm tip and the target point $\states^\target$, and $\chi_2 > 0$ is a suitably chosen regularization parameter. 

\begin{remark}
Careful analysis is needed regarding the controllability aspect of this infinite dimensional system. The Lie algebra rank condition or otherwise known as the Chow-Rashevsky theorem for finite dimensional systems  \cite{wei1939uber, rashevsky1938connecting, sussmann1972controllability} typically does not hold for infinite dimensional systems, and one needs additional assumptions, e.g. \cite{heintze1999homogeneity, salehani2014controllability}. 
Moreover, existence of the first order Pontryagin's Maximum Principle (PMP) type optimality conditions in the infinite dimensional settings is non-trivial. A few attempts have been made to show generalized PMP conditions for infinite dimensional systems with additional assumptions \cite{krastanov2011pontryagin, li2012optimal, halder2019optimality}. However, the scope of this paper is not to address these questions, rather to characterize optimal trajectories for a soft arm manipulation task, in a quest to explain experimentally observed behaviors. We will therefore proceed assuming that the controllability and PMP optimality conditions hold.
\end{remark}


\section{Optimal control solution} \label{sec:solution}
\subsection{The maximum principle} \label{sec:pmp}

The \textit{costate} is denoted as $\costates(t) := (\mu (t), \gamma (t)) \in \statespace^*, ~ t \in [0, T]$.  The control Hamiltonian function\footnote{Notice the difference between the Hamiltonian function $H$ in the optimal control theory and the Hamiltonian $\Hamiltonian$ in the elastic rod theory.} $H: \statespace \times \controlspace \times \R \times \statespace^* \rightarrow \R$ is defined as
\begin{align}
\begin{split}
H(z(t), u(t), \xi_0, \xi(t)) &:= \xi_0 \Lagrangian(z(t), u(t)) \\ 
									  &\quad	+ \inner{\xi(t)}{f(z(t), u(t))}
\end{split}
\label{eq:pre-Hmailtonian}
\end{align}
The Hamilton's equations in the infinite-dimensional settings are as follows:
\begin{proposition}[Maximum Principle \cite{krastanov2011pontryagin, halder2019optimality}]
Let $\bar{u} \in \mathfrak{U}$ be an optimal control for problem \eqref{eq:optimal_control_problem_free_end} and $\bar{z}(t)$ be the corresponding optimal trajectory. Then, there exists a pair $(\bar{\xi}_0, \bar{\xi}(t)) \in \R \times \statespace^*, ~ t \in [0, T]$, such that $(\bar{\xi}_0, \bar{\xi}) \not\equiv 0, \bar{\xi}_0 \leq 0$, $\bar{\xi}$ satisfies the differential equation
\begin{align}
\frac{\dif \bar{\xi}}{\dif t} (t) = - \left(\frac{\delta f}{\delta z} \right)^\dagger  (\bar{z}(t), \bar{u}(t)) \; \bar{\xi}(t) - \bar{\xi}_0 \frac{\delta \Lagrangian}{\delta z} (\bar{z}(t), \bar{u}(t))
\label{eq:adjoint}
\end{align}
where $(\cdot)^\dagger$ denotes the adjoint operator. The pointwise maximization of the pre-Hamiltonian holds, i.e. 
\begin{align}
H(\bar{z}(t), \bar{u}(t), \bar{\xi}_0, \bar{\xi}(t)) \geq H(\bar{z}(t), v, \bar{\xi}_0, \bar{\xi}(t)) 
\label{eq:hamiltonian_maximization}
\end{align} 
for all $v \in \controlspace$ and for all $t \in [0, T]$. Moreover, $\bar{z}$ and $\bar{\xi}$ satisfy Hamilton's canonical equations
\begin{align}
\begin{split}
\frac{\dif \bar{z}}{\dif t} (t) &= \frac{\delta H}{\delta \xi} (\bar{z}(t), \bar{u}(t), \bar{\xi}_0, \bar{\xi}(t)) \\
\frac{\dif \bar{\xi}}{\dif t} (t) &= - \frac{\delta H}{\delta z} (\bar{z}(t), \bar{u}(t), \bar{\xi}_0, \bar{\xi}(t)) 
\end{split}
\label{eq:hamiltons_equations}
\end{align}
Furthermore, the vector $\bar{\xi}(T)$ satisfies the transversality condition
\begin{align}
\bar{\xi}(T) = - \frac{\delta \Phi}{\delta z}(\bar{z}(T)) 
\label{eq:transversality_condition}
\end{align}
\end{proposition}
 \vspace*{15pt}
 In the remainder of this paper, we will restrict ourselves in studying only the normal extremals, i.e. where $\bar{\xi}_0 \neq 0$ and can be normalized to $-1$. The explicit form of the Hamilton's equations as a set of six (forward) PDEs and six (adjoint) PDEs appears in Appendix \ref{appdx:explicit}.  

\subsection{Computing optimal control -- the forward-backward algorithm} \label{sec:forward_backward}
\begin{figure*}[!ht]
	\centering
	\includegraphics[width=\textwidth, trim = {0pt 0pt 0pt 0pt}]{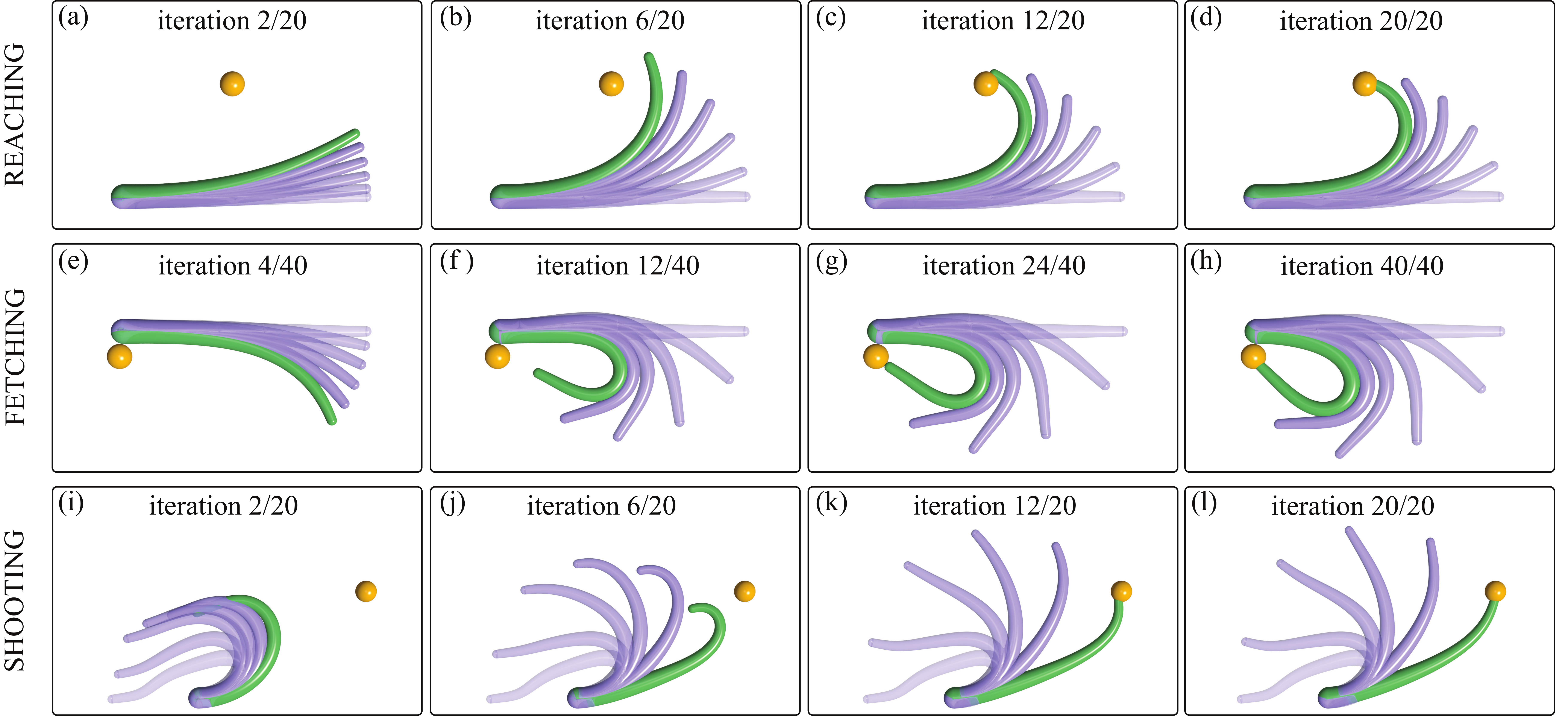}
	\caption{Summary of the numerical experiments: We select four iterations for each experiment. Six time instances, including the initial time $t=0$  and the terminal time $t=T$, are illustrated for each iteration. The rod at the terminal time is depicted in green while other time instances are depicted in fade-in purple. The target is represented by an orange ball. (a)-(d) The arm is initialized with straight, undeformed configuration and is tasked to reach the target located in the first quadrant at $r^\target=(9,9)$ [cm] with the tip. Simulation time is $T=0.5$ s for all 20 iterations. (e)-(h) The arm is initialized with straight, undeformed configuration and is tasked to reach the target located in the first quadrant at $r^\target=(0,-2)$ [cm] with the tip. Simulation time is $T=0.6$ s for all 40 iterations.  (i)-(l) The arm is initialized with bent, deformed configuration and is tasked to reach the target located at $r^\target=(16,10)$ [cm] with the tip. Simulation time is $T=0.8$ s for all 20 iterations.}
	\label{fig:cases}
	\vspace{-5pt}
\end{figure*}
A solution to the optimal control problem \eqref{eq:optimal_control_problem_free_end} necessarily has to satisfy the PMP conditions \eqref{eq:hamiltonian_maximization}, \eqref{eq:hamiltons_equations}, and \eqref{eq:transversality_condition}. This calls for solving the resulting two point boundary value problem in a function space. This is a challenging task even for a finite-dimensional nonlinear problem, for which various numerical techniques have been proposed~\cite{morrison1962multiple, bock1984multiple, von1993numerical}. 

An alternate approach is to employ an iterative algorithm (here referred to as forward-backward algorithm) to compute the optimal control. The idea is to start with an initial guess of the control $u^{(1)}$ in the first iteration.  (This guess may be zero.)  In each subsequent iteration, the control is modified so as to achieve the maximization of the control Hamiltonian $H$ \cite{bryson1962steepest, fujimoto2003iterative}. 

Suppose the state, costate and control at iteration $k$ is denoted as $z^{(k)}, \xi^{(k)}$, and $u^{(k)}$, respectively. At $k$-th iteration the steps of this algorithm are as follows:
\begin{enumerate}
\item
\textit{Run forward path:} The state equation \eqref{eq:control_system_abstract} is integrated forward in time from $t = 0$ to $T$, to obtain the state $z^{(k)}$. 
\item
\textit{Calculate terminal condition of the costate} from the transversality condition \eqref{eq:transversality_condition}.
\item
\textit{Run backward path:}
The costate, or the adjoint equation \eqref{eq:adjoint} is integrated backward in time from $t = T$ to $0$ to obtain the costate $\xi^{(k)}$.
\item
\textit{Update control:}
The triad $(z^{(k)}, \xi^{(k)}, u^{(k)})$ will typically not satisfy the Hamiltonian maximization criterion \eqref{eq:hamiltonian_maximization}.  Therefore, the control is updated in the direction of steepest ascent of the control Hamiltonian. Denoting the gradient of $H$ with respect to the control $u$ as $\frac{\delta H}{\delta u}$, the control update law is expressed as
\begin{align}
{u}^{(k+1)} &= {u}^{(k)} + \eta_k \frac{\delta H}{\delta {u}^{(k)}} 
\label{eq:forward_backward}
\end{align}  
where $\eta_k>0$ is the learning rate at iteration $k$. 

\end{enumerate}

Then we repeat steps 1)\,--\,4) until either of the two convergence criteria is met: i) the absolute change in control update becomes lower than a threshold $\epsilon$; ii) the number of iterations exceeds a predefined value.


%

\subsection{Numerical solver} \label{sec:solver}


Both the forward and backward path equations \eqref{eq:control_system_abstract}, \eqref{eq:adjoint} are systems of nonlinear PDEs that need to be propagated forward (or backward) in time given initial data. For the forward path, the specialized software~\textit{Elastica} \cite{gazzola2018forward} is used.  The software is designed for high-fidelity simulations of three dimensional Cosserat rods. A custom numerical solver is implemented for the backward adjoint equation. 


Both forward and backward dynamics solvers use finite difference techniques to discretize the spatial dimension. For the backward dynamics, certain spatial discretization operators are employed~\cite{bergou2008discrete,lang2011multi}, the details of which appear in the Appendix~\ref{appx:backward_numerics}. As for the time discretization, the forward dynamics are evolved via a position Verlet scheme. Such a scheme is commonly used to simulate a mechanical system where the state is decomposed into $(\states, \momentums)$ pair~\cite{verlet1967computer}. As explicit calculations show in Appendix \ref{appdx:explicit}, the costate $\xi$ is decomposed into a $(\mu, \gamma)$ pair which can be interpreted as velocity-position variables. Hence, the position Verlet scheme is also used for costate dynamics to integrate backward in time.

\section{Simulation results} \label{sec:numerics}
In this section, we demonstrate the numerical results of the optimal control on a single CyberOctopus arm of rest length $L_0$. In all our experiments, the intrinsic strains are chosen so that the arm is intrinsically straight, i.e. $\nu^\intrinsic = (1,0)$ and $\kappa^\intrinsic = 0$. The variable diameter $\phi(s) = \phi_{\text{base}} (L_0 - s) + \phi_{\text{tip}} s$ models the tapering of the arm. The cross sectional area and the second moment of area are given by $A = \tfrac{\pi \phi^2}{4}$ and $I = {A^2/}{4\pi}$. The effective shear modulus is given by $G=\frac{4}{3}\cdot\frac{E}{2(1+\text{Poisson's ratio})}$ \cite{gazzola2018forward}, where we take the Poisson's ratio to be 0.5 by assuming a perfectly incompressible isotropic material.
Parameters like density, modulus of elasticity, and physical dimensions are taken from \cite{yekutieli2005dynamic, chang2020energy}. Simulation parameters are tabulated in Table~\ref{tab:num_para}.

\begin{table}[t]
	\centering
	\caption{Parameters for Numerical Simulation}
	\begin{tabular}{ccc}
		\hline
		\hline\noalign{\smallskip}
		Parameter & Description & Numerical value \\
		\hline\noalign{\smallskip}
		\multicolumn{3}{c}{{\bf Rod model}}\\
		$L_0$ & length of the undeformed rod [cm] & $20$ \\
		$\phi_\text{base}$ & rod base diameter [cm] & $2$ \\
		$\phi_\text{tip}$ & rod tip diameter [cm] & $0.8$ \\
		$\rho$ & density [kg/${\text{m}}^3$] & $1042$ \\
		$\zeta$ & damping coefficient [kg/s]  & $0.01$ \\
		$E$ & Young's modulus [kPa] & $10$ \\
		\hline\noalign{\smallskip}
		\multicolumn{3}{c}{{\bf Numerics}}\\
		$\Delta t$ & Discrete time step-size [s] & $10^{-5}$ \\
		$N$ & number of discrete segments & $100$ \\
		$\epsilon$ & threshold for control convergence & $10^{-8}$ \\ 
		\hline 
	\end{tabular}
	\label{tab:num_para}
\end{table}

\subsection{Numerical experiments} \label{subsec:cases}
We test our solver to find the optimal trajectories for three different test cases.  We set the terminal tip orientation free and only penalize the distance between the terminal tip position and the target position $r^\target\in\R^2$, i.e. for $q^\target = (r^\target, \theta^\target)$, we use the following formula for $\Phi_{\text{tip}} (\cdot, \cdot)$ 
\begin{align}
\Phi_{\text{tip}} \left(q(T, L_0), q^\target \right) = \frac{1}{2} \norm{r(T, L_0) - r^\target}^2
\label{eq:phi_tip}
\end{align}
where the norm is the usual Euclidean distance in $\R^2$.

\subsubsection{Reaching task}
Our first experiment is a simple reaching problem. The arm is initialized to be straight and undeformed.  Our goal is to control the arm to reach the target with the tip at time $T = 0.5$ s. We consider the optimal control problem~\eqref{eq:optimal_control_problem_free_end}-\eqref{eq:terminal_cost},~\eqref{eq:phi_tip} with weight parameter $\chi_1=10$ and regularization parameter $\chi_2=2\times10^4$.  We ran the forward-backward algorithm for 20 iterations with fixed learning rate $\eta_k=3\times10^{-5}$. 

\begin{figure}[t]
	\centering
	\hspace{-10pt}
	\includegraphics[width=0.95\columnwidth, trim = {0 0 50 50}, clip = true]{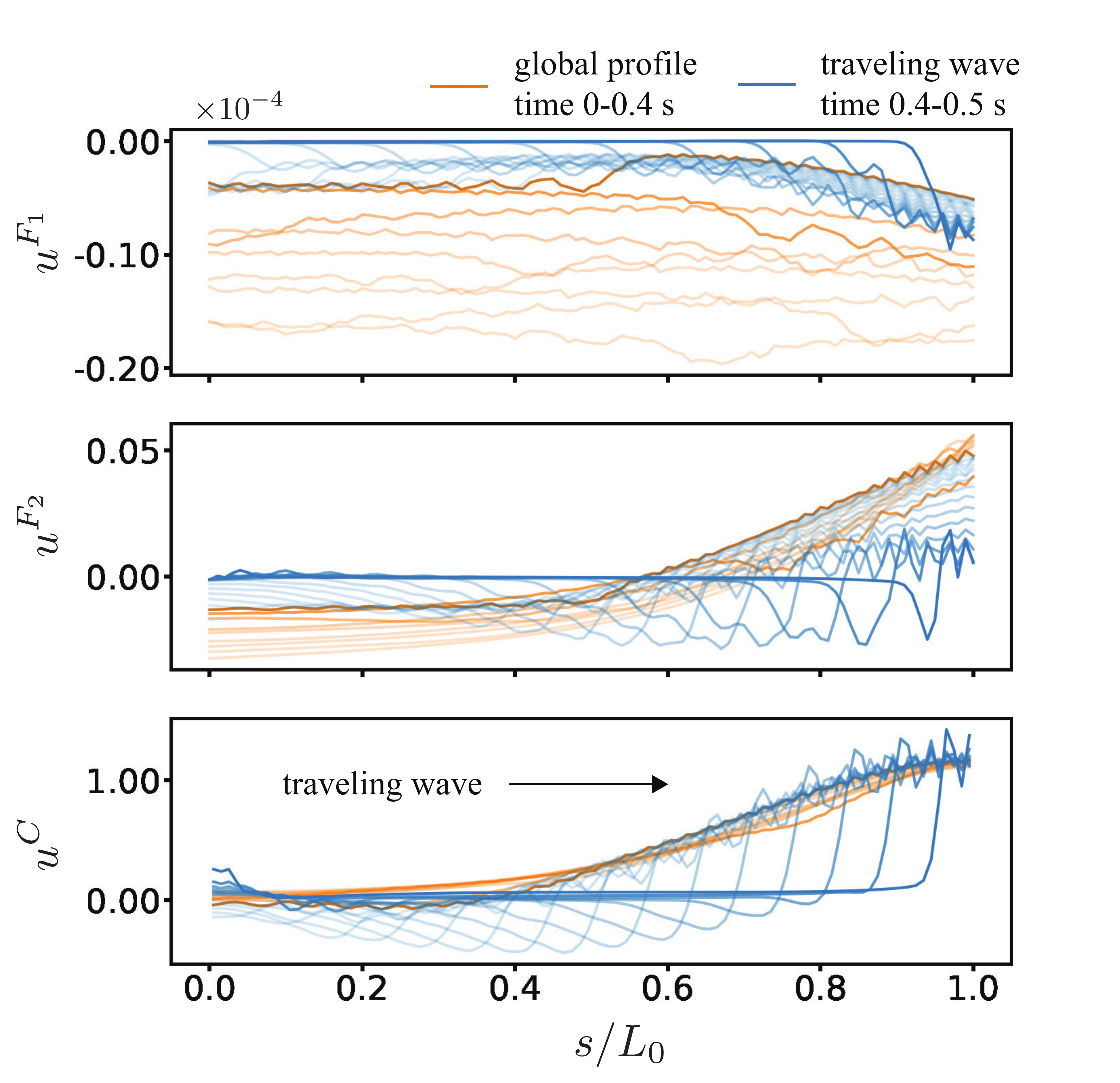} %
	\caption{Learned optimal controls for the reaching task: Control inputs $\uf = (u^{F_1}$, $u^{F_2})$ and $\uc$ along the arm are illustrated for the last iteration.  Nine time snapshots are shown in orange from $t=0$ s to $t=0.4$ s, and sixteen time snapshots are shown in blue from $t=0.4$ s to $t=0.5$ s (the most transparent lines correspond to the beginning of the time interval). 
The orange lines indicate the global profile of the optimal controls. The blue lines indicate the distinguishable traveling waves in optimal controls.}
	\label{fig:case1-ctrl} 
	\vspace{-10pt}
\end{figure}

We select four different iterations to demonstrate the control results. As we see in Fig~\ref{fig:cases}a-d, the reaching capability of the arm improves over iterations due to control updates. In the 2nd iteration, the arm {does not} bend much yet but shows the trend of moving towards the target. In the 6th iteration, the arm tip already gets close to the target. The controls converge quickly, and in the last iteration, the time snapshots show that the learned optimal control drives the arm to smoothly bend towards the target and the tip reaches the target at the terminal time. 
Fig.~\ref{fig:case1-ctrl} depicts the control inputs in the last iteration. We can see the emergence of a wave propagation in control inputs. 

\subsubsection{Fetching task}
During a fetching motion, the arm is observed to form several pseudo-joints ~\cite{sumbre2005motor}. To investigate this behavior, optimal trajectories are computed where the static target $r^\target$ is close to the base of the arm and is thought of as the mouth of the octopus. The arm is initialized to be straight and undeformed.  The forward-backward algorithm is run for 40 iterations with parameters $\chi_1=10$, $\chi_2=2\times10^4$ and $\eta_k=4\times10^{-5}$. The terminal time $T=0.6$ s is fixed for all iterations.
Fig.~\ref{fig:cases}e-h depicts the fetching movement where the arm forms a bend as it tries to get close to the target point. 


\subsubsection{Shooting task (reaching from bent position)}
Octopuses are known to curl up their arms while at rest, and when they try to catch food from a distance, they `shoot' one of the arms towards the target~\cite{sumbre2001control}. During this, the bend propagation is most prominently observed. 
Inspired by these observations, in our last experiment the arm is initialized at a bent position according to the initial curvature
\begin{equation*}
\kappa(0,s) = \sum_{i=1}^4 M_i\exp\left(-\frac{(s-s_i)^2}{2\times \sigma_i^2}\right)
\end{equation*} 
where $M_i$'s are [20, 78, 10, -30], $s_i$'s are [0, 0.3$L_0$, 0.7$L_0$, 0.85$L_0$], and $\sigma_i$'s are [0.015, 0.015, 0.012, 0.008]. 
Our goal is to reach the target at time $T = 0.8$ s. We ran the forward-backward algorithm for 20 iterations with parameters $\chi_1=100$, $\chi_2=2\times10^4$ and $\eta_k=3\times10^{-5}$. 

The control results of four forward-backward iterations are demonstrated in Fig.~\ref{fig:cases}i-l. Even though the arm reaches the target at the final iteration, the stereotypical bend propagation \cite{sumbre2001control} is absent. This is suggestive of the potential importance of environmental effects such as drag forces.

\subsection{Characteristics of optimal control} \label{subsec: tune-para}

{In our simulations, the optimal control solutions exhibit the following patterns. There is an initial global profile for the control. Starting from $t = 0$ s, a localized wave travels back and forth along the global profile and the magnitude of this wave increases as $t$ increases. At first, the wave is not discernible due to its small magnitude and thus, the global profile is dominant as indicated by the orange lines in Fig.~\ref{fig:case1-ctrl}. As time $t$ nears the final time $T$, the wave traveling from the base to the tip of the arm becomes more visible and it dominates the control as shown by the blue lines in Fig.~\ref{fig:case1-ctrl}.  

We vary different parameters to investigate how they affect the optimal control solution, especially the wave propagation. }

\subsubsection{Wave speed}
We observed that the parameters of the optimal control problem, e.g. $T, \chi_1, \chi_2, r^{\text{target}}$ (see discussion in Sec.~\ref{sec:chi_1} about the parameter $\chi_1$), geometry of the arm (the length and tapered diameter profile), dissipation constant $\zeta$, and numerical integration constants (e.g. $N, \Delta t, \eta$) do not affect the speed $c$ of the wave. However, Young's modulus ($E$) and density ($\rho$) of the arm do affect the wave speed. We calculate the speeds for different sets of $E$ and $\rho$ values, which shows the linear relationship (the graphic is omitted due to lack of space) 
\begin{equation*}
	c = 0.653\sqrt{\tfrac{E}{\rho}}
	\label{eq:wave_speed}
\end{equation*}
This experiment indicates that the wave in the optimal control solution is actually a fundamental property of the elastic arm. Further study is required to draw connections to the stereotypical bend propagation waves observed in octopuses \cite{sumbre2001control}.

\subsubsection{Parameter $\chi_1$}
 \label{sec:chi_1}
\begin{figure}[!t]
	\centering
	\hspace{-20pt}
	\includegraphics[width=0.95\columnwidth, trim = {10 10 90 60},  clip = true]{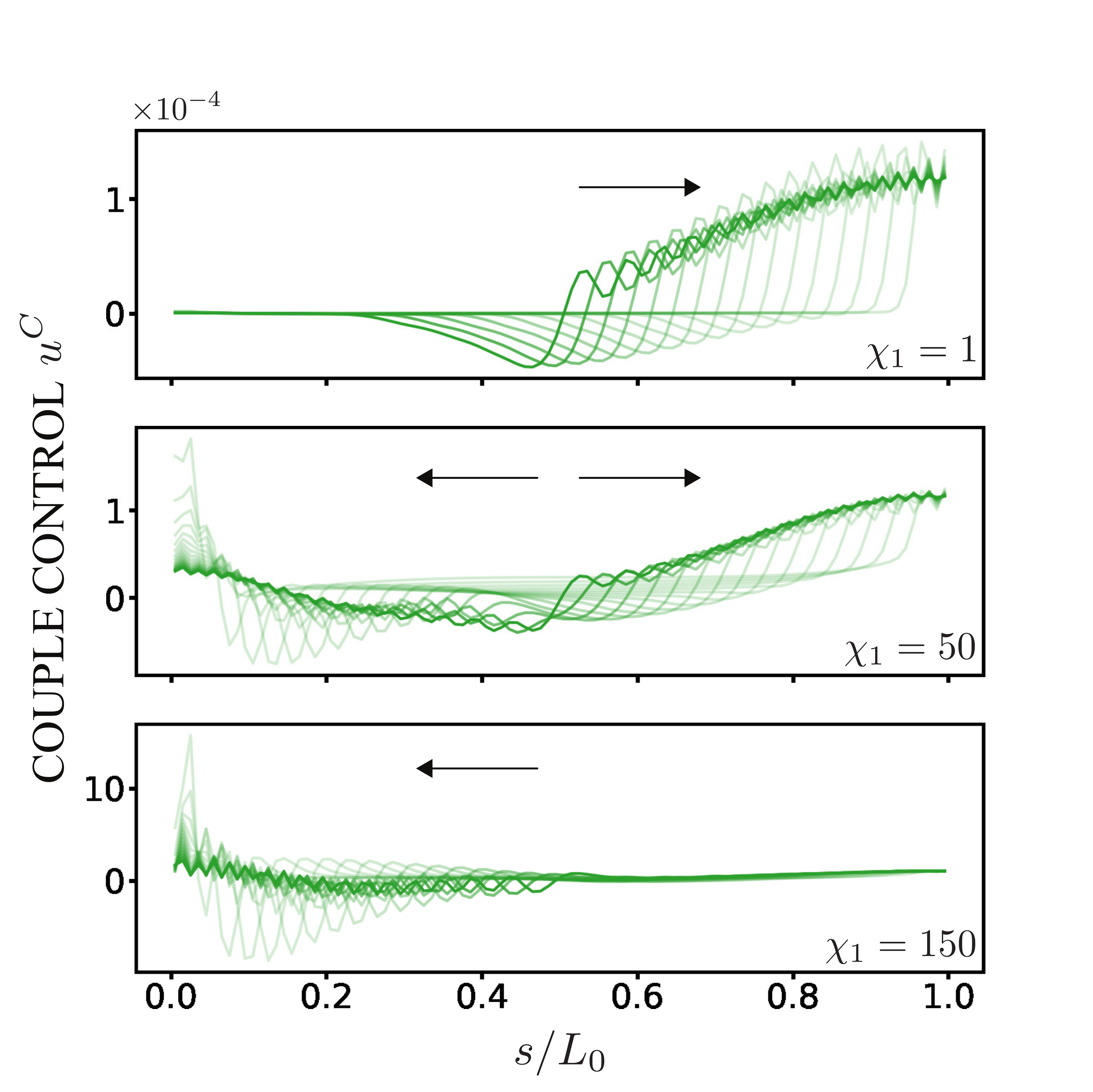} 
	\caption{Comparison of wave behaviors in couple control for different $\chi_1$ parameters in the reaching task: Sixteen time snapshots are shown in green from $t=0.45$ s (most solid) to $t=0.5$ s. (most transparent) The black arrows indicate the direction of the dominant wave propagation. For $\chi_1=1$, the usual dominating wave travels from base to the tip.  For $\chi_1=50$, both the original wave and a second wave are visible and they meet at the middle point $0.5L_0$. For $\chi_1=150$, the second wave is dominant which travels from tip to the base.}
	\label{fig:two_waves}
	\vspace{-10pt}
\end{figure}

In the cost function~\eqref{eq:Langrangian}, we penalize the deformation of the arm with the parameter $\chi_1$. Even though the parameter $\chi_1$ does not affect the wave speed, increasing $\chi_1$ leads to an interesting observation. When $\chi_1$ is high enough, a visible second wave appears in the control solution which propagates in the opposite direction of the original wave. Moreover, these two waves meet exactly at the middle point of the arm (Fig.~\ref{fig:two_waves}). The resemblance of this behavior with the observation of \cite{sumbre2005motor} demands further analysis.

\section{Conclusion and Future Work} \label{sec:conclusion}

In this paper, we investigate an optimal control problem for a single CyberOctopus arm modeled as a planar Cosserat rod. A free endpoint optimal control problem is formulated to minimize the control energy and a weighted potential energy of the rod. To reach a target point, the proximity of the arm's tip to the target point is penalized at the terminal time. The necessary first order optimality conditions yield two systems, the Cosserat rod dynamics (forward) and the adjoint dynamics (backward), both described by nonlinear PDEs. To numerically solve these PDEs, specific spatial and temporal discretization techniques are used. The optimal controls are found by updating the controls in an iterative manner called the forward-backward algorithm. 
This framework is used to solve several biologically motivated control tasks. These numerical experiments reveal emergence of propagating waves in the optimal controls. However, the stereotypical bend propagation along the arm is not discovered under our current problem formulation. This motivates us to consider environmental effects like drag, and constraints of muscle actuation into our optimal control framework.

\bibliographystyle{IEEEtran}
\bibliography{reference}

\appendices
\renewcommand{\thelemma}{A-\arabic{section}.\arabic{lemma}}
\renewcommand{\thetheorem}{A-\arabic{section}.\arabic{theorem}}
\renewcommand{\theequation}{A-\arabic{equation}}
\renewcommand{\thedefinition}{A-\arabic{definition}}
\setcounter{lemma}{0}
\setcounter{theorem}{0}
\setcounter{equation}{0}

\section{Explicit calculations} \label{appdx:explicit}
\subsection{Details of a planar Cosserat rod dynamics}
For the planar case of the Cosserat rod, we denote $q=(r,\theta)$ as the state where the position vector along the rod $r(t,s) \in \R^2$ and the angle $\theta(t, s) \in\R$ can be used to measure local strains -- stretch $(\nu_1)$, shear $(\nu_2)$, and curvature $(\kappa)$. These are defined as follows:
\begin{equation*}
	r_s = Q\nu,\quad \theta_s = \kappa
\end{equation*}
where $Q =\begin{psmallmatrix*}[c]	\cos\theta & -\sin\theta \\ \sin\theta & \cos\theta
\end{psmallmatrix*}$ is the planar rotation matrix, and $\nu = \begin{psmallmatrix}\nu_1 \\ \nu_2 \end{psmallmatrix}$. The internal stresses, i.e. the forces $n$ (represented in the {material frame}) and couple $m$ are related to the stored energy function $W$ by
\begin{align*}
n = \frac{\partial W}{\partial \nu}, \quad m = \frac{\partial W}{\partial \kappa}
\end{align*}
We take the following quadratic form of $W$ so that the stress-strain relationship becomes linear
\begin{align*}
W = \frac{1}{2}\left( (\nu - \nu^\intrinsic)^{\mathsf{T}}S (\nu - \nu^\intrinsic) + B(\kappa - \kappa^\intrinsic)^2 \right)
\end{align*}
where the intrinsic strains of the rod are denoted by $(\nu^\intrinsic, \kappa^\intrinsic)$. Here, $S=\diag(EA,GA)$ is the stretch-shear rigidity matrix and $B = EI$ is the bending rigidity. $E, G$ are the Young's modulus and shear modulus, respectively.
 
Let us denote $\pt=\rho Ar_t$ and $\pr=\rho I\theta_t$ as the momentum variables $\momentums = (\pt, \pr)$, where $\rho$ is the density, $A$ is the cross sectional area and $I$ is the second moment of area. Let `$\cdot$' denote the dot product of two planar vectors, and `$\times$' represent the component of the cross product of two planar vectors along the normal vector that is coming out of the plane, i.e. $\begin{psmallmatrix}	x_1 \\ x_2 \end{psmallmatrix} \cdot 
\begin{psmallmatrix}	y_1 \\ y_2 \end{psmallmatrix} = x_1y_1 + x_2y_2$ , and $\begin{psmallmatrix}	x_1 \\ x_2 \end{psmallmatrix} \times
\begin{psmallmatrix}	y_1 \\ y_2 \end{psmallmatrix} = x_1y_2 - x_2y_1$.

The Cosserat dynamics \eqref{eq:control_system_abstract} are written as 
\begin{equation}
	\begin{aligned}
		r_t &= \frac{1}{\rho A}\pt \\
		\theta_t &= \frac{1}{\rho I}\pr \\
		\pt_t &= (Q n)_s - \frac{1}{\rho A}\zeta\pt + \uf \\
		\pr_t &= (m)_s + \nu \times n - \frac{1}{\rho I}\zeta\pr + \uc
	\end{aligned}
\end{equation}
where $u = (\uf, \uc)$ denote the force and couple control inputs. 
%
%

\subsection{Details of the adjoint equations}

Denote the costate to $(q,p) = ((r,\theta),(\pt,\pr))$ as $(\mu, \gamma) = ((\mut,\mur),(\gammat,\gammar))$. Then, the pre-Hamiltonian \eqref{eq:pre-Hmailtonian} is explicitly written as
\begin{equation}
\footnotesize
	\begin{aligned}
		H &= \int_0^{L_0} \left[ \frac{1}{\rho A}\mut\cdot\pt + \frac{1}{\rho I}\mur\pr + \gammat\cdot\left((Q n)_s - \frac{1}{\rho A}\zeta\pt\right) \right.\\
		&\quad \left.+ \gammar\left((m)_s + \nu\times n - \frac{1}{\rho I}\zeta\pr\right) \right. + \gammat\cdot\uf + \gammar\uc\\
		&\quad \left. - \frac{1}{2}\left(\uf\cdot\uf + \left(\uc\right)^2\right) - \chi_1 \mathcal{V}(q) \right]\ud s 
	\end{aligned}
\normalsize
\end{equation}
Maximizing $H$ with respect to $u$ gives the first order necessary condition for optimal control
\begin{align}
\uf = \gammat,\quad \uc = \gammar
\label{eq:optimal_u_cosserat}
\end{align}
Furthermore, the costate evolution equations \eqref{eq:adjoint} take the explicit form
\begin{equation}
\footnotesize
	\begin{aligned}
		&\begin{aligned}
			\mut_t &= -\frac{\delta H}{\delta r} \\ 
			&= -\left(QSQ^\mathsf{T}\gammat_s\right)_s - \left[QM_1(GA\nu - EA\sigma)\gammar\right]_s - \chi_1(Qn)_s
		\end{aligned} \\
		&\begin{aligned}
			\mur_t = &-\frac{\delta H}{\delta \theta} \\
			= &-\left(B\gammar_s\right)_s + \left[Q\left(M_2 n - SM_2\nu\right)\right]\cdot  \gammat_s \\ & + \left[(M_2\nu)\times n+\nu\times(SM_2\nu)\right]\gammar - \chi_1\left((m)_s + \nu\times n\right)
		\end{aligned} \\
		&\begin{aligned}
			\gammat_t &= -\frac{\delta H}{\delta \pt} = -\frac{1}{\rho A} \left(\mut - \zeta\gammat\right)
		\end{aligned} \\
		&\begin{aligned}
		\gammar_t &= -\frac{\delta H}{\delta \pr} = -\frac{1}{\rho I} \left(\mur - \zeta\gammar\right)
		\end{aligned}
	\end{aligned}
\normalsize
	\label{eq:costate_cosserat}
\end{equation}
where $\sigma = \nu - \nu^\intrinsic$, $M_1=\begin{psmallmatrix}
	0 & 1 \\ 1 & 0
\end{psmallmatrix}$ and $M_2=\begin{psmallmatrix}
	0 & -1 \\ 1 & 0
\end{psmallmatrix}$.

These equations are to be accompanied with the transversality condition \eqref{eq:transversality_condition} (with \eqref{eq:terminal_cost}, \eqref{eq:phi_tip})
\begin{equation}
	\begin{aligned}
		\mut(T,s) &=  - \delta(s-L_0)\left[\chi_2(r(t, s) -r^\target)\right]\Big|_{t=T} \\
		\mur(T,s) &= 0 \\
		\gammat(T,s) &= 0 \\
		\gammar(T,s) &= 0
	\end{aligned}
	\label{eq:transversality_cosserat}
\end{equation}
where $\delta(\cdot)$ denotes the delta function.

\subsection{Control update law}
Denoting $u = (u^F, u^C)$ and $\gamma = (\gamma^r, \gamma^\theta)$, we can write the control update law~\eqref{eq:forward_backward} for the forward-backward algorithm at iteration $k$ as
\begin{equation}
\begin{aligned}
{u}^{(k+1)} &= {u}^{(k)} + \eta_k \frac{\delta H}{\delta{u}^{(k)}} = {u}^{(k)} + \eta_k\left({\gamma}^{(k)} - {u}^{(k)} \right)
\end{aligned}
\label{eq:control_update}
\end{equation}

\section{Numerical Methods} \label{appx:backward_numerics}

We use the following spatial and temporal discretization for the backward path that is consistent with the forward path.

\subsection{Spatial discretization}

In the software package \textit{Elastica}, the Cosserat rod is decomposed into $N+1$ nodes for the position $r$ and $N$ segments for the angle $\theta$~\cite{gazzola2018forward}. 

We define the following two difference operators for vectors according to finite difference approximation~\cite{bergou2008discrete,lang2011multi}. Let $\{\R^p\}_N$ denote a set of $N$ vectors in $\R^p$. Then, 
$\D:\{\R^p\}_N\mapsto\{\R^p\}_{N+1}$ and $\diff:\{\R^p\}_N\mapsto\{\R^p\}_{N-1}$ are defined as follows:
\begin{equation}
\footnotesize
\y_{i=1,\ldots,N+1} = \D(\x_{j=1,\ldots,N}) = \left\{\begin{aligned}
& \x_1,\qquad\qquad\qquad\quad \ i=1 \\
& \x_i - \x_{i-1},\quad i=2,\ldots,N \\
& -\x_N,\quad\quad\quad\ i=N+1
\end{aligned}\right.
\normalsize
\end{equation}
and 
\begin{equation}
\footnotesize
\begin{aligned}
c_{\ell=1,\ldots,N-1} &= \diff(\x_{j=1,\ldots,N}) = \x_{\ell+1} - \x_{\ell},\quad \ell=1,\ldots,N-1
\end{aligned}
\normalsize
\end{equation}
where $\y_i\in\R^p$ for $i=1,\ldots,N+1$, $\x_j\in\R^p$ for $j=1,\ldots,N$ and $c_\ell\in\R^q$ for $\ell=1,\ldots,N-1$. Note that $\D$ and $\diff$ operate on a set of $N$ vectors and then return $N+1$ and $N-1$ vectors, respectively.

Now for the rest of this Appendix, we will use specific subscripts $(\cdot)_i$, $(\cdot)_j$ and $(\cdot)_\ell$ to denote the set of discretized variables with the dimension of spatial discretization to be $N+1$, $N$ and $N-1$, respectively.

For the backward path, we discretize the costate into $\mut_i$, $\gammat_i$ and $\mur_j$, $\gammar_j$.
Then the first-order necessary condition for optimal control is
\begin{equation}
\begin{aligned}
\uf_i &= \gammat_i \\
\uc_j &= \gammar_j
\end{aligned}
\label{eq:optimal_u_cosserat-discrete}
\end{equation}
where $\uf_i$ and $\uc_j$ are the discretized control inputs to be used in the forward path.

The costate dynamics~\eqref{eq:costate_cosserat} are discretized as follows:
\begin{equation}
\footnotesize
\begin{aligned}
\frac{\ud \mut_i}{\ud  t} = &-\D\left(Q_jSQ_j^\mathsf{T}\diff(\gammat_i)/\Delta s\right) - \D\Big(Q_jM_1(GA\nu_j - EA\sigma_j)\gammar\Big) \\ &- \chi_1 \D\left(Q_j n_j\right) \\
\frac{\ud \mur_j}{\ud  t} = &-\D\left(B\diff(\gammar_j)/\Delta s\right) + \left[Q_j\left(M_2 n_j - SM_2\nu_j\right)\right] \cdot \diff(\gammat_i) \\ &+ \left[(M_2\nu_j)\times n_j+\nu_j\times(SM_2\nu_j)\right]\gammar_j \Delta s \\ &- \chi_1 \left(\D\left(m_\ell\right) + (\nu_j\times n_j)\Delta s\right) \\
\frac{\ud \gammat_i}{\ud  t} = &-\frac{1}{\rho A} \left(\mut_i - \zeta\gammat_i\right) \\
\frac{\ud \gammar_j}{\ud  t} = &-\frac{1}{\rho I} \left(\mur_j - \zeta\gammar_j\right) 
\end{aligned}
\normalsize
\label{eq:costate_cosserat-discrete}
\end{equation}
where $\Delta s=L_0/N$ is the length of each discretized segment of the rod. $r_i$ , $Q_j$, $\nu_j$, $\sigma_j$, $n_j$ and $m_\ell$ are discretized variables obtained from the forward path. Details of these variables are covered in~\cite{gazzola2018forward}.

The transversality conditions~\eqref{eq:transversality_cosserat} are discretized into
\begin{equation}
\begin{aligned}
\mut_i(T) &= - \delta(i-(N+1))\left[\chi_2(r_i-r^\target)\right]\Big|_{t=T} \\
\mur_j(T) &= 0 \\
\gammat_i(T) &= 0 \\
\gammar_j(T) &= 0
\end{aligned}
\label{eq:transversality_cosserat-discrete}
\end{equation} 


\subsection{Time discretization}

We use the second-order position Verlet time integration~\cite{gazzola2018forward} as follows:
\begin{equation}
\begin{aligned}
\gammat_i\left(t-\frac{\Delta t}{2} \right) &= \gammat_i(t) - \frac{\Delta t}{2} \frac{\ud \gammat_i}{\ud t}(t) \\
\mut_i(t-\Delta t) &= \mut_i(t) - \Delta t \frac{\ud \mut_i}{\ud t}\left(t-\frac{\Delta t}{2} \right) \\
\gammat_i(t-\Delta t) &= \gammat_i\left(t-\frac{\Delta t}{2}\right) - \frac{\Delta t}{2} \frac{\ud \gammat_i}{\ud t}(t-\Delta t)
\end{aligned}
\label{eq:time-space-discrete}
\end{equation}
Similarly for $\gammar_j$ and $\mur_j$.

\end{document}